\documentclass{article}

\usepackage{url}

\usepackage{amsmath, amssymb, amstext, amsopn, amsxtra, graphicx, color, calligra, hyperref, url, ulem, polynom, permute, ifthen, float, epstopdf, linsys, marvosym, rotating}

\begin{document}

\newcommand{\recip}[1]{\ensuremath{\frac1{#1}}}
\newcommand{\half}{\recip2}
\newcommand{\third}{\recip{3}}
\newcommand{\twothirds}{\ensuremath{\frac23}}
\newcommand{\fourth}{\recip4}
\newcommand{\fifth}{\recip5}
\newcommand{\sixth}{\recip6}
\newcommand{\point}[1]{\mbox{\(\left(#1\right)\)}}
\newcommand{\qar}{\begin{eqnarray*}}
\newcommand{\raq}{\end{eqnarray*}}
\newcommand{\inv}{^{-1}}
\newcommand{\vone}{{\bf 1}}
\newcommand{\va}{{\bf a}}
\newcommand{\vb}{{\bf b}}
\newcommand{\vr}{{\bf r}}
\newcommand{\vs}{{\bf s}}
\newcommand{\vu}{{\bf u}}
\newcommand{\comment}[1]{\hskip 5mm \parbox{3in}{#1}}
\newcommand{\casI}{\renewcommand{\labelenumi}{\bf Case \Roman{enumi}.}\num}
\newcommand{\sacI}{\mun \renewcommand{\labelenumi}{\arabic{enumi}.}}
\newcommand{\num}{\begin{enumerate}}
\newcommand{\mun}{\end{enumerate}}
\newcommand{\rowr}{\ensuremath{\left(\begin{array}{rrrr}r_1 & r_2 & \cdots & r_k \end{array}\right)}}
\newcommand{\rows}{\ensuremath{\left(\begin{array}{rrrr}s_1 & s_2 & \cdots & s_k \end{array}\right)}}
\newcommand{\qed}{\mbox{\(\square\,\,\)}}
\newcommand{\longcomment}[1]{}

\newcommand{\nakaut}{Nakayama automorphism}
\newcommand{\scg}{\ensuremath{\dk/\dk^2}}
\newcommand{\Rd}{\ensuremath{\hat R}}
\newcommand{\Hom}[1]{\mbox{\rm Hom}_{#1}}
\newcommand{\leftmods}{\mbox{left $R$-modules}}
\newcommand{\bimod}{$(R,R)$-bimodule}
\newcommand{\bimods}{$(R,R)$-bimodules}
\newcommand{\Aut}{\mbox{\rm Aut}}
\newcommand{\End}{\mbox{\rm End}}
\newcommand{\dk}{\ensuremath{\dot{k}}}
\newcommand{\m}{\ensuremath{\mathfrak m}}
\newcommand{\isom}{\simeq}
\newcommand{\isommap}{\stackrel{\sim}{\longrightarrow}}
\newcommand{\lomar}{\cite{lomar}}
\newcommand{\Rrad}{\ensuremath{R/\mbox{rad}R}}
\newcommand{\R}{{_R}}
\newcommand{\lRd}{\ensuremath{\R\hat R}}
\newcommand{\hi}{\ensuremath{\varphi}}
\newcommand{\sig}{\ensuremath{\sigma}}
\newcommand{\alp}{\ensuremath{\alpha}}
\newcommand{\el}{\ensuremath{\ell}}
\newcommand{\ro}{\ensuremath{\rho}}
\newcommand{\rad}{\ensuremath{\mbox{rad }R}}
\newcommand{\all}{ \hskip 5mm (\forall}
\newcommand{\Id}{\mbox{\rm Id}}
\newcommand{\field}[1]{\ensuremath{\mathbb #1}}
\newcommand{\CC}{\field C}

\newtheorem{thm}{Theorem}
\newtheorem{cor}[thm]{Corollary}
\newtheorem{conj}[thm]{Conjecture}
\newtheorem{lem}[thm]{Lemma}
\newtheorem{rem}[thm]{Remark}
\newtheorem{exmp}[thm]{Example}

\renewcommand{\em}{\textit}

\title{Bilinear Forms on Frobenius Algebras}

\author{Will Murray\\
wmurray@csulb.edu\\
California State University, Long Beach\\
Long Beach, CA 98040-1001}
\date{January 24, 2014}

\maketitle

\begin{abstract}
We analyze the homothety types of associative bilinear forms that can occur on a Hopf algebra or on a local Frobenius \(k\)-algebra \(R\) with residue field \(k\).  If \(R\) is symmetric, then there exists a unique form on \(R\) up to homothety iff \(R\) is commutative.  If \(R\) is Frobenius, then we introduce a norm based on the \nakaut\ of \(R\).  We show that if two forms on \(R\) are homothetic, then the norm of the unit separating them is central, and we conjecture the converse.  We show that if the dimension of \(R\) is even, then the determinant of a form on \(R\), taken in \scg, is an invariant for \(R\).  

\end{abstract}

\section{Introduction}

Let $R$ be a finite-dimensional algebra over a field $k$.  We say $R$ is a \em{Frobenius algebra} if there exists a nondegenerate bilinear form \(B: R \times R \rightarrow k\) that is \em{associative} in the sense that \(B(rs,t) = B(r,st) \; \forall r,s,t \in R\).  We say $R$ is a \em{symmetric algebra} if there exists a nondegenerate associative symmetric bilinear form \(B: R \times R \rightarrow k\).  These properties are equivalent to the existence of an isomorphism between $R$ and its $k$-dual $\Rd := \Hom{k}(R,k)$ as \leftmods, respectively, as \bimods.  


Since many different isomorphisms between \(R\) and \(\Rd\) exist, we also have many bilinear forms.  A natural question to ask is whether the various forms are \em{isometric}, that is, the same under change of basis.  It is trivial to observe that any form may be scaled by a nonzero constant from \(k\), so we define two forms \(B\) and \(B'\) to be \em{homothetic} if there exists a change of basis \(V \in \Aut_k(R)\) and a scalar \(\alpha \in \dk := k -\{0\}\) such that \(B'(r,s) = \alpha B(Vr, Vs) \; \forall r,s \in R\).  We then ask instead when two forms on \(R\) are homothetic.

In this paper we will study the question above in the case when \(R\) has an ideal \m\ with \(R/\m \isom k\).  For example, this condition is satisfied by the group algebra \(R = kG\), where \(G\) is any finite group, by taking \m\ to be the kernel of the augmentation map \(\epsilon(\sum \alpha_g g) = \sum \alpha_g \in k\).  It is also true for Hopf algebras, whose definition includes the existence of a counit \(\epsilon: R \rightarrow k\).  For most of our results we will also need to assume that \(k\) has good characteristic.  

We will show that in the local symmetric case, there exists a unique form on \(R\) up to homothety iff \(R\) is commutative.  For Frobenius algebras that are not symmetric, we will introduce a norm based on the order of the \nakaut, a distinguished $k$-algebra automorphism of \(R\) that measures how far $R$ is from being a symmetric algebra.  (The automorphism is the identity iff $R$ is symmetric.)  We will show that if two forms on \(R\) are homothetic, then the norm of the unit separating them is central, and we will conjecture the converse.  Finally, we will study the determinant of a form on \(R\) and show that in even dimension, the value of the determinant in \scg\ is an invariant of the algebra.  

The idea of comparing an algebra with its dual was pioneered by F. G. Frobenius himself (\cite{frobenius}) in connection with representations of finite groups, and group algebras have remained important examples of symmetric algebras.  Nakayama gave new examples and developed the main properties of Frobenius algebras and symmetric algebras in \cite{nak1,nak2,nak3}.  More recently, the group algebra example was generalized when Larson and Sweedler (\cite{larson}) showed that all finite dimensional Hopf algebras are Frobenius.  (See \cite{bergman} for a treatment of the ubiquity of Hopf algebras.)  Modern interest in Frobenius algebras has grown far beyond their representation-theoretic origins as connections have been discovered to such diverse areas as topological quantum field theories, Gorenstein rings in commutative algebra, coding theory, and the Yang-Baxter Equation.  For an excellent reference on the subject, see \lomar.  

\section{Preliminaries and examples}\label{section-preliminaries}

Let \(k\) be a field and \(R\) be a finite-dimensional \(k\)-algebra.  Throughout this paper, we will use the word \em{form} (respectively, \em{symmetric form}) to mean a nondegenerate bilinear form (respectively, nondegenerate symmetric bilinear form) \(B: R \times R \rightarrow k\) that is \em{associative} in the sense that \(B(rs,t) = B(r,st) \; \forall r,s,t \in R\).

In \lomar, Theorems 3.15 and 16.54, we have:

\begin{thm} \label{theorem-frobdef}The following conditions are equivalent:
\begin{description}
\item[1.]  \(R \isom \Rd\) as \leftmods\ (respectively, as \bimods).
\item[2.]  There exists a linear functional \(\lambda: R \rightarrow k\) whose kernel contains no nonzero 
left ideals.  (Respectively, \(\lambda(rs) = \lambda(sr) \; \forall r,s \in R\).)
\item[3.]  There exists a form (respectively, symmetric form) \(B: R \times R \rightarrow k\).  
\end{description}
If \(R\) satisfies these conditions, \(R\) is said to be a Frobenius (respectively, symmetric) algebra.\qed
\end{thm} 

For any finite-dimensional algebra \(R\), the dual \Rd\ is isomorphic as a left $R$-module to the injective hull of \Rrad, so the isomorphism type of \lRd\ does not depend on the ground field \(k\).  In (\cite{murray}) it is shown that the isomorphism type of \Rd\ as an \bimod\ is also independent of \(k\); thus we may speak of \(R\) being a Frobenius (respectively, symmetric) algebra without reference to the ground field.  

The equivalence of the first two conditions in Theorem~\ref{theorem-frobdef} follows from taking \(\lambda\) to be the image of 1 under the module isomorphism and vice versa.  The equivalence of the last two follows from defining \(B(r,s):=\lambda(rs)\) and \(\lambda(r):= B(r,1)\).  Since the third condition is right-left symmetric, we could also include the right-handed analogues of the other conditions above.  

Given one isomorphism \(\hi:\)\(_R R \isommap \lRd\), any other isomorphism $\hi'$ is obtained by composition with an automorphism of the left regular module $_R R$, which corresponds to right multiplication by a unit $u \in U(R)$.  This affects the other conditions above as follows:  the new functional is \(\lambda' = u\lambda: r \mapsto \lambda(ru)\); and the new form is \(B'(r,s) = B(r,su)\).

If the conditions of Theorem~\ref{theorem-frobdef} hold, the nondegeneracy of the form $B$ implies that there is a unique $k$-linear map \(\sigma: R \rightarrow R\) defined by \(B(r,s)=B(s,\sigma(r)) \; \forall r,s \in R\).  (Equivalently, \(\lambda(rs) = \lambda(s\sig(r)) \; \forall r,s \in R\).)  It is easy to check that $\sigma$ is actually a $k$-algebra automorphism of $R$, known as the \em{\nakaut} of $R$.  Replacing $B$ with a new form $B'$ defined by the unit $u$ gives us the new automorphism \(\sigma' = I_u \circ \sig\), where \(I_u\) is the inner automorphism \(r \mapsto uru\inv\).  So the \nakaut\ is determined up to composition with inner automorphisms; equivalently, it is a well-defined element of the group of outer automorphisms of $R$.  The algebra is symmetric iff \sig\ can be taken to be the identity, iff the \nakaut\ determined by an arbitrary form is an inner automorphism.

In most theorems we will assume that \(R\) has an ideal \m\ with \(R/\m \isom k\).  As shown in the previous section, this condition is satisified when \(R\) is a finite-dimensional Hopf algebra or when \(R\) is the group algebra of any finite group over any field.  

We will sometimes assume additionally that \(R\) is local, or, equivalently, that \(\m = \rad\).  This condition is satisfied, for example, by the group algebra \(kG\) where \(G\) is a finite \(p\)-group and char \(k = p\).

Another classic example that satisfies both assumptions above is the Nakayama-Nesbitt example.  This algebra is described in \cite{nak1} in terms of matrices, but we will use polynomials to save space.  Let \alp\ be a fixed element of \dk.  We define
\[
R_\alpha := k\langle x,y\rangle /(x^2 = y^2 = 0, yx = \alp xy) = k \oplus kx \oplus ky \oplus kxy,
\]
a four-dimensional local algebra with maximal ideal \(\m = (x,y)\).  The functional \(\lambda(a + bx + cy + dxy) = d\) shows that \(R_\alpha\) is Frobenius.  The \nakaut\ is \(\sig: x \mapsto \alp \inv x, y \mapsto \alp y\), since, for example, \(\lambda(x\sig(y)) = \lambda(yx) = \lambda(x (\alp y))\).  The algebra is symmetric iff \(\alp = 1\), and the \nakaut\ has finite order iff \alp\ is a root of unity in \(k\).  

In Example~\ref{example-xyx} we will extend this example by increasing the index of nilpotency of \m.  It is also possible to increase the number of generators of \m, in which case the relations can be arranged so that the \nakaut\ effects any desired permutation among the generators.  

\section{The transpose of an endomorphism}

Let \(R\) be a Frobenius algebra with form \(B\).  Given a $k$-vector space endomorphism \(\hi \in \End(R_k)\), the nondegeneracy of \(B\) gives us a unique endomorphism \(\hi^{t}\), the {\em transpose} of \hi\ with respect to \(B\), that satisfies
\[
B(\hi r,s) = B(r, \hi^{t} s) \all r,s \in R).
\]

\longcomment{
If we fix a $k$-basis for \(R\), then we can abuse notation slightly and think of elements of \(R\) as column vectors and \(B\) as a square matrix, so that
\[
B(r,s) = r^T B s,
\]
where \(T\) represents ordinary transpose of matrices.  Again abusing notation and thinking of \hi\ as a square matrix with respect to the same basis, we obtain
\[
\hi^{t} = B\inv \hi^T B.
\]
} 

We have the standard properties:
\begin{eqnarray*}
(\hi+\psi)^{t} & =&  \hi^{t} + \psi^{t}, \\
(\hi\psi)^{t} & = &  \psi^{t} \hi^{t}, \\
(\alpha \hi)^{t} & = & \alpha (\hi^{t}) \all \alpha \in k), \\
(\hi\inv)^{t} & = & (\hi^{t})\inv \mbox{ if $\hi$ is invertible.}
\end{eqnarray*}
In general, however, it is not true that \((\hi^{t})^{t} = \hi\), unless \(B\) is symmetric.  We will see that it depends on the order of the \nakaut\ \sig.  We will use the notation \(\hi^{t^2}:= 
(\hi^{t})^{t}\), and so on.

\begin{lem}\label{lemma-t2}
\(\hi^{t^2} = \sigma \hi \sigma\inv.\)
\end{lem}
\textit{Proof}.  
For all $r,s \in R,$
\[
B(r,\hi^{t^2} s) = B(\hi^{t} r, s) = B(\sigma\inv s, \hi^{t} r) = B(\hi \sigma\inv s, r) = B(r, \sigma \hi \sigma\inv s).  \qed
\]

\begin{cor}\label{cor-transposecycle}
If \(\sigma^n = \Id,\) then \(\hi^{t^{2n}} = \hi\).  \qed
\end{cor}

Now for any \(x \in R\), let \(\el_x, \ro_x \in \End_k(R)\) be left and right multiplication by \(x\), respectively.

\begin{lem}\label{lemma-transpose}
For all \(i \geq 0\), 
\begin{itemize}
\item \(\ro_x^{t^{2i}} = \ro_{\sigma^i(x)}\).
\item \(\ro_x^{t^{2i+1}} = \el_{\sigma^i(x)}\).
\end{itemize}
\end{lem}

\textit{Proof}.  
For all $r,s \in R,$
\begin{eqnarray*}
B(r, \ro_x^{t} s) & = & B(\ro_x r, s) \\
	& = & B(rx,s) \\
	& = & B(r,xs) \\
	& = & B(r,\el_x s), 
\end{eqnarray*}
proving the second statement for \(i=0\).  Similarly, for all \(r,s \in R\),
\begin{eqnarray*}
B(r, \ro_x^{t^2} s) & = & B(\ro_x^{t} r, s) \\
	& = & B(xr,s) \mbox{ by the above}\\
	& = & B(x,r s)  \\
	& = & B(r s, \sigma(x)) \\
	& = & B(r, s \sigma(x)) \\
	& = & B(r, \ro_{\sig(x)} s),
\end{eqnarray*}
proving the first statement for $i=1$.  The general statements follow by induction.  \qed

\begin{cor}\label{lemma-transposecommute}
Even ``transpose powers'' of \(\ro_x\) commute with odd 
``transpose powers'' of \(\ro_x\).  \qed
\end{cor}

Now let \(B'\) be another form on \(R\) such that \(B'(r,s) = B(r,su) \; \forall r,s \in R\), where \(u \in U(R)\).  The transpose gives us a way to determine when the two forms are homothetic:  

\begin{lem}\label{lemma-vtv}
\(B\) and \(B'\) are homothetic iff \(\exists \alpha \in \dk\) and \(V \in \Aut_k(R)\) such that \mbox{\(\ro_u = \alp V^{t} V \in \Aut_k(R)\)}.
\end{lem}

\textit{Proof}.  
If \(\ro_u = \alp V^{t} V\), then for all \(r, s \in R\), 
\[
B'(r,s) = B(r,su) = B(r,\ro_u s) = B(r, \alp V^{t} V s) = \alp B(Vr, Vs),
\]
so \(B\) and \(B'\) are homothetic.  The converse is similar.  \qed

\section{Forms on local symmetric algebras}

Throughout this section we will assume that char \(k \neq 2\) and that \(R\) is a local symmetric algebra with maximal ideal \m\ such that \(R/\m \isom k\).  This condition is satisfied, for example, by the group algebra \(kG\) where \(G\) is a finite \(p\)-group and char \(k = p > 2\), or by the Nakayama-Nesbitt example above.  We will show that up to homothety, \(R\) has a unique symmetric form, and it has a unique form iff it is commutative.

Some of the results here are similar to those derived by Watanabe in \cite{watanabe}.  Watanabe studied the three-fold multilinear form \(\theta\) on a symmetric algebra defined by \(\theta (r,s,t):= \lambda(rst)\), where \(\lambda\) is as in Theorem~\ref{theorem-frobdef}.  He showed (by elementary techniques) that algebras with isometric three-fold forms are isomorphic.  (This is done without using our assumptions that the algebra be local or that char \(k \neq 2\).)  Watanabe also studied bilinear forms and is responsible for Lemma~\ref{lemma-watanabe} below.

Throughout this section, let \(B\) and \(B'\) be two forms on $R$ such that \(B'(r,s) = B(r,su) \; \forall r,s \in R\), where \(u \in U(R)\).  As above, let \(\ro_u, \el_u \in \Aut_k(R)\) be right and left multiplication by \(u\).  

\begin{lem}\label{lemma-watanabe}
(Lemma 2 in \cite{watanabe})  If \(B\) is symmetric, then \(B'\) is symmetric iff \(u \in Z(R)\).  
\end{lem}

\textit{Proof}.  
For all \(r,s \in R\), we have
\begin{eqnarray*}
B'(r,s) = B'(s,r) & \iff & B(r,su) = B(s,ru)\\
& \iff & B(r,su) = B(ru,s)\\
& \iff & B(r,su) = B(r,us)\\
& \iff & B(r,\ro_u s) = B(r,\el_u s).
\end{eqnarray*}
Using the nondegeneracy of \(B\), we have that \(B'\) is symmetric iff \(\ro_u = \el_u\), iff \(u \in Z(R)\).  \qed

\begin{lem}\label{lemma-symisom}
Suppose char \(k \neq 2\) and $(R,\m)$ is a local symmetric $k$-algebra with $R/\m \isom k$.  If \(B\) is symmetric, then \(B\) and \(B'\) are homothetic iff \(u \in Z(R)\).
\end{lem}

\textit{Proof}.  
If \(B\) and \(B'\) are homothetic, then \(B'\) is symmetric, so by Lemma~\ref{lemma-watanabe}, \(u \in Z(R)\).

Conversely, suppose that \(u \in Z(R)\); we claim that we can find \(\alp \in \dk\) and \(v \in Z(R)\) such that \(\alp v^2 = u\).  If we can find such a \(v\) (necessarily in \(U(R)\) since \(u \in U(R)\)), then we have
\[
\ro_u = \alp (\ro_v)^2 = \alp \el_v \ro_v = \alp \ro_v^{t} \ro_v,
\]
so \(B\) and \(B'\) are homothetic by Lemma~\ref{lemma-vtv}.

To prove the claim, suppose that \(u = \alpha + m \in R\) where \(\alpha \in \dk\) and \(m \in \m\).  We set \(v:= 1 + a_1 m + a_2 m^2 + \cdots,\) where \(a_i \in k\), and solve inductively for the \(a_i\) to satisfy
\begin{eqnarray*}
\alp v^2 & = & u\\
\alp + 2 \alp a_1 m + \alp(2a_2 + a_1^2)m^2 + \cdots & = & \alp + m.
\end{eqnarray*}
(Note that \(m^n = 0\) for some \(n\), so the expression above will terminate.)  So it suffices to solve 
\begin{eqnarray*}
2 \alp a_1 & = & 1 \\
(2  a_2 + a_1^2) & = & 0 \\
(2  a_3 + 2a_1 a_2) & = & 0 \\
& \vdots & .
\end{eqnarray*}
Since \(2 \alp \in \dk\), these equations can always be solved, proving the claim.  \qed

It is now easy to prove our main results on symmetric algebras.

\begin{thm}\label{theorem-uniquesym}
Suppose char \(k \neq 2\) and $(R,\m)$ is a local symmetric $k$-algebra with $R/\m \isom k$.  Then there exists a unique symmetric form on $R$ up to homothety.  
\end{thm}

\textit{Proof}.  
Suppose \(B\) and \(B'\) are both symmetric forms on $R$.  Then by Lemma~\ref{lemma-watanabe}, $u \in Z(R)$, so by Lemma~\ref{lemma-symisom}, \(B\) and \(B'\) are homothetic.  \qed

\longcomment{
\begin{remark}
Although we will not present the proof here, as it would take us too far away from our goal of analyzing forms on local algebras, Theorem~\ref{theorem-uniquesym} extends to the nonlocal case in which \Rbar is a matrix ring over $k$.
\end{remark}
}

\begin{thm}\label{theorem-uniquecom}
Suppose char $k \neq 2$ and $(R,\m)$ is a local symmetric $k$-algebra with $R/\m \isom k$.  Then there exists a unique homothety class of forms on $R$ iff $R$ is commutative.  
\end{thm}

\textit{Proof}.  
From Lemma~\ref{lemma-symisom}, the forms \(B\) and \(B'\) are homothetic iff $u \in Z(R)$.  

If $R$ is commutative, then certainly \(u \in Z(R) \; \forall u \in U(R)\), so all forms are homothetic.  (Alternatively, we could prove this implication by noting that any form on a commutative algebra is symmetric, so Theorem~\ref{theorem-uniquesym} shows that all forms are homothetic.)  

Conversely, if all forms are homothetic, then \(u \in Z(R) \; \forall u \in U(R).\)  This implies that $R$ is commutative, since (using the fact that $R$ is local) for any $r \in R$, one of $r$ or $1+r$ is a unit, hence central, so $r$ is central.  \qed

We remark that Frobenius algebras that are not symmetric may have a unique homothety class of forms without being commutative.  In fact, the noncommutative Nakayama-Nesbitt algebra above does have a unique homothety class of forms.  So we cannot omit the symmetric condition from the theorem above.  

Theorem~\ref{theorem-uniquecom} is a special case of Theorem~\ref{theorem-centralnorm} and 
Conjecture~\ref{conjecture-converse}, which treat the case in which the \nakaut\ has finite order instead of being the identity.  

\section{Forms on Frobenius algebras}

We will try to generalize the results of the previous section to Frobenius algebras.  We would like to focus on the case in which the \nakaut\ of $R$ has finite order.  However, since the \nakaut\ is only well-defined up to inner automorphism, we first examine the case when it has finite inner order, i.e.\ it has finite order as a member of the group of outer automorphisms, i.e.\ some power of any particular \nakaut\ is inner.  In this case we will define a norm on $R$ and use it to show that in good characteristic we can then find a new form whose \nakaut\ actually has finite order.  Then we will show that a necessary (and perhaps sufficient) condition for any other form to be homothetic to that form is that the norm of the associated unit be central. 

Suppose $R$ is a Frobenius $k$-algebra with form $B$ and corresponding \nakaut\ \sig.  Recall from Section~\ref{section-preliminaries} that if $B'(r,s) = B(r,su)$, then $B'$ has \nakaut\ \(\sigma' = I_u \circ \sigma\), where $I_u$ denotes the inner automorphism $r \mapsto uru\inv.\)  Hence, the inner order of the \nakaut\ is independent of the choice of the form.  We use this to define our norm function on $R$:  suppose $\sigma$ has finite inner order $n$.  Then for $r \in R$, we define \(N_\sigma(r):= r\sigma(r)\cdots\sigma^{n-1}(r).\)  Before using this norm to find another \nakaut\ that has finite order, we need some technical lemmas.  

\begin{lem}\label{lemma-fixsigma}
Suppose \(\sigma^n = I_a.\)  Then \(\sigma(a) = a.\) 
\end{lem}

\textit{Proof}.  
For all \(r\in R,\)
\begin{eqnarray*}
B(r,a) & = & B(a, \sigma(r)) = B(\sigma(r), \sigma(a)) = \cdots\\
\cdots & = & B(\sigma^n(r), \sigma^n(a)) = B(ara\inv, a) = B(ar, 1) = B(a,r).  \qed
\end{eqnarray*}

\begin{lem} \label{lemma-innerorder}
Suppose \(\sigma^n = I_a.\)  Then \((\sigma')^n = I_{N_\sigma(u)a}.\)
\end{lem}

\textit{Proof}.  
Let $r \in R.$  Then
\begin{eqnarray*}
\sigma'(r) & = & I_u \circ \sigma(r) = u \sigma(r) u\inv \\
(\sigma')^2(r) & = & u \sigma(u) \sigma^2(r) \sigma(u\inv) u\inv \\
	& \vdots & \\
(\sigma')^n(r) & = & u \sigma(u) \cdots \sigma^{n-1}(u) \sigma^n(r) \sigma^{n-1}(u\inv)\cdots \sigma(u 
\inv) u \inv \\
	& = & N_\sigma(u) I_a(r) (N_\sigma(u))\inv \\
	& = & I_{N_\sigma(u)} \circ I_a (r) \\
	& = & I_{N_\sigma(u)a}(r).  \qed
\end{eqnarray*}

\begin{lem}\label{lemma-isnorm}
Suppose $(R,\m)$ is local Frobenius with $R/\m \isom k$.  Suppose \(\sigma^n = I_a\) and char $k \nmid \,n$.  If $x \in U(R)$ satisfies \(\sigma(x) = x\) and \(\bar{x} \in \dk^n\) in $R/\m$, then there exists \(u \in U(R)\) such that \(\sigma(u) = u\) and $x = u^n$, so \(N_\sigma(u) = x\).
\end{lem}

\textit{Proof}.  
Let \(x = \alp + m\) with \(\alp \in \dk, m \in \m,\) and say \(\alp = u_0^n\) for some \(u_0 \in \dk.\)  We set \(u = u_0 + u_1 m + u_2 m^2 + \cdots\) (this terminates since $m$ is nilpotent) and solve inductively for \(u_i \in k\) to satisfy \(u^n = x\), that is,
\[
u_0^n + nu_0^{n-1}u_1 m + \left[n u_0^{n-1} u_2 + {n \choose 2} u_0^{n-2} u_1^2 \right] m^2 + \cdots  =  \alp + m.
\]
It suffices to solve
\begin{eqnarray*}
u_0^n & = & \alpha \\
nu_0^{n-1}u_1 & = & 1 \\
n u_0^{n-1} u_2 + {n \choose 2} u_0^{n-2} u_1^2 & = & 0 \\
& \vdots & .
\end{eqnarray*}
Since $n \neq 0$ in $k$, these equations can be solved for the $u_i$.  Then \(u^n = x\) and \(u \in U(R)\) since \(u_0 \neq 0\).  Finally, since \(\sigma(x) = x\) and \(\sigma = \Id\) on $\alpha \in k$, we have \(\sigma(m) = m\), so \(\sigma(u) = u\) and hence \(N_\sigma(u) = u^n = x\).  \qed

\begin{thm}\label{theorem-innerform}
Suppose $(R,\m)$ is local Frobenius with $R/\m \isom k$.  Suppose \(\sigma^n = I_a\) and char $k \nmid \,n$.  Then there exists a form $B'$ on $R$ whose \nakaut\ $\sigma'$ satisfies \((\sigma')^n = \Id.\)
\end{thm}

\textit{Proof}.  
Say \(a = \alpha + m\) with \(\alpha \in \dk, m \in \m.\)  We may assume that \(\alp = 1\), since scalar multiples of \(a\) will yield the same inner automorphism.  
Now by Lemma~\ref{lemma-fixsigma}, \(\sigma(a) = a,\) so \(\sigma(a\inv) = a\inv\).  Since \(a\inv = 1 - m + m^2 - \cdots\) and \(1 \in k^n,\) Lemma~\ref{lemma-isnorm} tells us that \(a\inv = N_\sigma(u)\) for some \(u \in U(R).\)  Let $B'$ be the form defined by $B'(r,s) = B(r,su)$ for \(r,s \in R\); then by Lemma~\ref{lemma-innerorder},
\[
(\sigma')^n = I_{N_\sigma(u) a} = I_{a\inv a} = \Id.  \qed
\]

Thus, in good characteristic we can find a form whose 
automorphism has finite order.  

We remark in passing that in the Ore ring of right twisted polynomials \(R[t,\sig]\) (where \(tr = \sig(r)t \;\forall r \in R\)), evaluation of the monomial \(t^i\) at \(r \in R\) is defined as \(N_i(r) := r \sig(r) \cdots \sig^{i-1}(r)\).  (See \cite{lamleroy}.)  Thus our norm \(N_\sigma(r)\) coincides with evaluation of \(t^n\) at \(r\).  


\longcomment{  
In fact, the theorem above has a direct analogue for Ore polynomials:

\begin{thm}\label{theorem-twisted}
Let $R$ be a ring with an automorphism $\sigma$ and suppose that \(\sigma^n = I_a\) for some $a \in U(R).$  If the 
twisted polynomial \(t^n - a\inv \in R[t,\sigma]\) has a root in $R$, then there exists an automorphism $\sigma'$ on 
$R$ such that \(R[t,\sigma] \isom R[t',\sigma']\) and \((\sigma')^n = \Id\) (so in particular, \((t')^n \in 
Z(R[t',\sigma'])\)).
\end{thm}

\textit{Proof}.  
Let $u \in R$ be the root, so \(N_n(u) = N_\sigma(u) = a\inv.\)  Setting \(\sigma':= I_u 
\circ \sigma,\) we have \((\sigma')^n = \Id\) as in Theorem~\ref{theorem-innerform}, and then
\(R[t',\sigma'] \isom R[t,\sigma]\) via \(t' \mapsto ut.\)  \qed

Analogous results on differential polynomial rings have been proven by T. Y. Lam in \cite{charp}.
}

Returning to Frobenius algebras, we have established that in many cases we can find a form whose \nakaut\ has finite order.  We now ask when other forms will be homothetic to that form.  We will show that a necessary (and perhaps sufficient) condition for this to occur is that the unit separating the two forms have a central norm.  

As before, let $B$ and $B'$ be two forms on $R$ such that \(B'(r,s) = B(r,su)\), where $u \in U(R)$.  Let $B$ have Nakayama automorphism $\sigma$, and suppose \sig\ has finite order $n$.  We saw in Lemma~\ref{lemma-symisom} that if $n=1$ (so $R$ is symmetric, and $N_\sigma(u) = u$) then $B'$ is homothetic to $B$ iff $u \in Z(R)$.  We now generalize that result.  We can relax the condition that \(R\) be local and assume only that \(R/\m \isom k\).  This assumption is satisfied by the group algebra of any finite group over any field and by any finite-dimensional Hopf algebra.  

\begin{thm}\label{theorem-centralnorm}
Suppose $R$ is Frobenius with an ideal \m\ such that $R/\m \isom k$, and suppose $\sigma^n = \Id.$  If $B$ and $B'$ are homothetic, then \(N_\sigma(u) \in Z(R).\)
\end{thm}

\begin{conj}\label{conjecture-converse}
The converse is true, too, i.e.\ if $R$ is Frobenius with an ideal \m\ such that $R/\m \isom k$ and $\sigma^n = \Id,$ then \(N_\sigma(u) \in Z(R)\) implies that $B$ and $B'$ are homothetic.
\end{conj}

Experimental evidence supports this conjecture, at least in the local case.  

\textit{Proof of Theorem~\ref{theorem-centralnorm}}.  
Let $\ro_u \in \Aut_k(R)$ be right multiplication by $u$.  If $B$ and $B'$ are homothetic, then by Lemma~\ref{lemma-vtv}, \(\exists \alpha \in \dk,  V \in \Aut_k(R)\) such that \mbox{\(\ro_u = \alp V^{t} V  \in \Aut_k(R)\)}.  Then for \(i \geq 0\), we have
\begin{eqnarray*}
\ro_u^{t^{2i}} & = & \alp V^{t^{2i+1}} V^{t^{2i}}\\
\ro_u^{t^{2i+1}} & = & \alp V^{t^{2i+1}} V^{t^{2i+2}}.
\end{eqnarray*}

We now expand the expression
\begin{equation}  \label{equation-manyros}
\ro_u\inv \ro_u^{t} \left(\ro_u^{t^2}\right) \inv \ro_u^{t^3} \dots \left(\ro_u^{t^{2n-2}}\right) \inv \ro_u^{t^{2n-1}},
\end{equation}
noting that all the \(\alp\)'s cancel and that \(\ro_u^{t^{2n-1}} = \alp V^{t^{2n-1}} V^{t^{2n}} = \alp V^{t^{2n-1}} V\) by Corollary~\ref{cor-transposecycle}.  Expression~\eqref{equation-manyros} then becomes
\[
V\inv \left(V^t\right) \inv \cdot V^t V^{t^2} \cdot \left(V^{t^2}\right)\inv \left(V^{t^3}\right)\inv \cdot \cdots \cdot 
 \left(V^{t^{2n-2}}\right)\inv \left(V^{t^{2n-1}}\right)\inv \cdot V^{t^{2n-1}} V = \Id.
\]
Now by Lemma~\ref{lemma-transposecommute}, the even transpose powers appearing in Expression~\eqref{equation-manyros} commute past the odd ones, and we can use Lemma~\ref{lemma-transpose} to translate these back into multiplication in $R$.
\begin{eqnarray*}
\left[\ro_u\inv \left(\ro_u^{t^2}\right) \inv \cdots \left(\ro_u^{t^{2n-2}}\right) \inv\right]\left[ \ro_u^{t} \ro_u^{t^3} \cdots \ro_u^{t^{2n-1}} \right] & = & \Id\\
\ro_u^{t} \ro_u^{t^3} \cdots \ro_u^{t^{2n-1}} & = & \ro_u^{t^{2n-2}} \cdots \ro_u^{t^2} \ro_u\\
\ro_u^{t} \ro_u^{t^3} \cdots \ro_u^{t^{2n-1}} r &  = & \ro_u^{t^{2n-2}} \cdots \ro_u^{t^2} \ro_u r \all r \in R)\\
u \sigma(u) \cdots \sigma^{n-1}(u) r & = & r u \sigma(u) \cdots \sigma^{n-1}(u) \all r \in R)\\
N_\sigma(u) r & = & r N_\sigma(u) \all r \in R)
\end{eqnarray*}
so \(N_\sigma(u) \in Z(R)\), as desired.  \qed

To illustrate the conjecture of the converse, we offer an algebra with two forms that are not homothetic for which the norm of the unit separating them is not central.

\begin{exmp}\label{example-xyx}
Consider the extended Nakayama-Nesbitt algebra
\[
R := \CC\langle x,y\rangle /(x^2 = y^2 = 0, \, xyx = yxy), 
\] 
a six-dimensional algebra with $\CC$-basis \(\{1, x, y, xy, yx, xyx\}\).  Then $R$ possesses a pair of forms that are not mutually homothetic, separated by a unit whose norm is not central.
\end{exmp}

\textit{Proof}.  
As with the original Nakayma-Nesbitt example, the functional \linebreak \mbox{\(\lambda (a + bx + cy + dxy + eyx + 
fxyx):= f\)} and the corresponding form \linebreak  \mbox{\(B(r,s) = \lambda(rs)\)} show that $R$ is Frobenius.  The \nakaut\ is  \mbox{\(\sigma: x \mapsto y, y \mapsto x\)}, because, for example $B(x, yx) = 1 = B(yx, y)$, so $y = \sigma(x).$  Then, of course, $\sigma^2 = \Id$.

Fix $\epsilon \in \CC$ and let \(u := 1 + \epsilon x \in U(R)\).  Then we have a form $B_\epsilon$ defined by \(B_\epsilon(r,s):= B(r,su)\).  Fixing the ordered basis \(\{1, x, y, xy, yx, xyx\}\), we abuse the notation slightly and think of elements of $R$ as column vectors and \(B_\epsilon\) as a square matrix.  Then \(B_\epsilon (r,s) = r^T B_\epsilon s\), where $T$ denotes ordinary matrix transposition.  The matrix for $B_\epsilon$ is
\[ 
B_\epsilon = \left( \begin{array}{cccccc}
0 & 0 & 0 & \epsilon & 0 & 1  \\
0 & 0 & \epsilon & 0 & 1 & 0  \\
0 & 0 & 0 & 1 & 0 & 0  \\
\epsilon & 1 & 0 & 0 & 0 & 0  \\
0 & 0 & 1 & 0 & 0 & 0  \\
1 & 0 & 0 & 0 & 0 & 0  \end{array} \right).
\]
For example, there is an \(\epsilon\) in the $(2,3)$ position because 

\[
B_\epsilon(x,y) = B(x,y(1 + \epsilon x)) = \lambda(xy(1+ \epsilon x)) = \lambda(xy + \epsilon xyx) = \epsilon.
\]

We can use the matrix to calculate the new \nakaut\ in matrix form, which we will denote by \(\Sigma_\epsilon\):
\[
s^T B_\epsilon \Sigma_\epsilon r = r^T B_\epsilon s = s^T B_\epsilon^T r\all r,s \in R),
\]
so \(\Sigma_\epsilon = B_\epsilon\inv B_\epsilon^T\).   

We claim, and we will prove below, that the forms for \(\epsilon = 0\) and \(\epsilon \neq 0\) are not homothetic.  Conjecture~\ref{conjecture-converse} then predicts that \(N_\sigma(u) \not\in Z(R)\) for \(\epsilon \neq 0\), and indeed,
\[
N_\sigma(u) = u \sigma(u) = (1+ \epsilon x) (1 + \epsilon y) = 1 + \epsilon x + \epsilon y + \epsilon^2 xy \not \in Z(R)
\]
because, for example, it does not commute with $x$.  

To prove the claim, we note that if a form \(B'\) is homothetic to \(B\), then the matrix for its \nakaut\ must be similar to the matrix for the automorphism for \(B\).   This is because if \(B'(r,s) = \alp B(Vr, Vs) \; \forall r,s \in R\), then \(B' = \alp V^T B V\), so 
\[
(B')\inv (B')^T = V\inv B\inv (V^T)\inv V^T B^T V = V\inv B\inv B^T V.
\]
However, it may be checked (by comparing Jordan Canonical Forms) that the matrix \(B_\epsilon\inv B_\epsilon^T\) when \(\epsilon \neq 0\) is \em{not} similar to the matrix \(B_\epsilon\inv B_\epsilon^T\) when \(\epsilon = 0\); therefore the forms are not homothetic.  This confirms Conjecture~\ref{conjecture-converse} for this particular pair of forms.  \qed

We note that although the unit above does not have a central norm, the algebra does contain nontrivial units with central norms.  If we take $u = 1 + \epsilon x - \epsilon y$ for some $\epsilon \in \dot\CC$, then 
\[
N_\sigma(u) = (1 + \epsilon x - \epsilon y)(1 - \epsilon x + \epsilon y) = 1 + \epsilon^2 xy + \epsilon^2 yx \in 
Z(R).
\]
If Conjecture~\ref{conjecture-converse} is true, then the form $B'$ defined by $u$ should be homothetic to $B$.  And indeed, this turns out to be true, although we suppress the somewhat laborious calculations necessary to confirm this.  
\section{Determinants of forms}

Although the results above show that the homothety class of a form on an algebra may not be well-defined, we will show now that for even-dimensional local Frobenius algebras, the determinant of a form is an invariant for the algebra.  

We first note that if \(B\) is a form on the Frobenius $k$-algebra \(R\), then \(\det B\) is a well-defined element of the square class group \scg.  To see this, fix an ordered $k$-basis for \(R\) and, by abuse of notation, write \(B\) for the matrix of the form with respect to this basis.  If we use the map \(V \in \Aut_k R\) to change bases, then the form \(B(Vr, Vs)\) will have matrix \(V^T B V\), which has the same determinant as \(B\) in \scg.  

Now suppose \(B'\) is another form on \(R\) separated from \(B\) by the unit \(u \in U(R)\), so \(B'(r,s) = B(r,su) \; \forall r,s \in R\).  Writing \(\ro_u\) for right multiplication by \(u\), the new form has matrix \(B' = B\ro_u\), so \(\det B' = \det B \det \ro_u\).  

If \(u = \alp \in \dk\), then \(\ro_u = \alp \Id\), which has determinant $\alpha^m$, where \(m = \dim_k R\).  If $m$ is odd, then \(\overline{\alp^m} = \bar{\alp} \in \scg\), so \(\det B'\) could be anything and it is pointless to hope for any significance to the determinant.  In the even-dimensional case, though, it turns out to be an invariant for the algebra:

\begin{thm}\label{theorem-determinant}
Let $k$ be a field of arbitrary characteristic and let $(R,\m)$ be an even-dimensional local Frobenius $k$-algebra 
with $R/\m \isom k$ and form $B$.  Then \(\det B\), valued in  \scg, is independent of the choice of the form.  Given another even-dimensional local Frobenius algebra \((R', \m'\)) with \(R'/\m' \isom k\) and form \(B'\), we have \(\det B = \det B'\) in \scg\ if \(R\) and \(R'\) are isomorphic as $k$-algebras.
\end{thm}
\textit{Proof}.  
By the discussion above it suffices to show that $\det \ro_u$ is a square in $\dk$.  We will do this essentially by showing that the matrix for $\ro_u$ with respect to a suitable basis is upper triangular.  Suppose \(\m^n \neq 0 = \m^{n+1}\).  Construct an ordered basis \(\{e_i\}\) for $R$ as follows:  start with a basis for \(\m^n\); then complete it to a basis for \(\m^{n-1}\), and so on.  

Let \(u = \alp (1 + m)\), where \(\alpha \in \dk, m \in \m\).  Then we may assume that \(\alp = 1\), since multiplying \(\ro_u\) by a scalar matrix will only change its determinant by a square.  

In the matrix of $\ro_u $ with respect to the basis $\{e_i\}$, the $i$-th column is given by the coordinates of $\ro_u e_i$ with respect to the basis.  Now
\[ 
\ro_u e_i = e_i u = e_i (1+m) = e_i + e_i m.
\]
Since \(e_i m\) lies in a strictly higher power of \m\ than \(e_i\), it is in the span of basis vectors \(e_j\) with \(j < i\).  Hence the $i$-th column of $\ro_u $ has a 1 in the $i$-th row, 0's below the $i$-th row, and other undetermined entries above.  So $\ro_u $ is upper triangular with 1's on the diagonal.  So $\det \ro_u  = 1$, as desired.  \qed

\begin{exmp}\label{example-determinant}
Let $a, b, c \in \dk$ with $b^2 \neq ac$, and let $R$ be the four-dimensional local commutative 
$k$-algebra
\[
R := k[x,y]/(ax^2 = bxy = cy^2, (x,y)^3 = 0) = k \oplus kx \oplus ky \oplus kx^2.
\]
Then the determinant of any form on $R$, up to square, is \(\delta(R) = ac(ac - b^2)\).  
\end{exmp}

Thus, if two such algebras have different values for \(\delta\) (taken in the square class group 
\(\dk/\dk^2\)) we know immediately that they are not isomorphic.  

\textit{Proof}.  
We have a form \(B(r,s) = \lambda(rs)\), where \(\lambda (a + bx + cy + dx^2) := d\).  The matrix for $B$ with respect to the ordered basis \(\{1, x, y, x^2\}\) is
\[
\left( \begin{array}{llll} 0 & 0 & 0 & 1 \\
0 & 1 & \frac{a}{b} & 0 \\
0 & \frac{a}{b} & \frac{a}{c} & 0 \\
1 & 0 & 0 & 0 \end{array} \right).
\]
The determinant is $a^2/b^2 - a/c$, which in the square class group is equal to $ac(ac - b^2)$.  \qed

In this example, the converse is true too:  the value of \(\delta\) in \scg\ actually determines the algebra up to isomorphism.  To see this, consider the algebra 
\[
R' := k[u,v]/(u^2 = - \delta v^2, uv = 0, (u,v)^3 = 0) = k \oplus ku \oplus kv \oplus ku^2,
\]
which is isomorphic to \(R\) under the map \(u \mapsto \delta y, v \mapsto abx - acy\).  It is then clear that the isomorphism type of \(R'\) is unchanged if the value of \(\delta\) is multiplied by a nonzero square.

\small
\noindent {\bf Acknowledgement}  The author would like to thank T.\ Y.\ Lam for useful conversations about the material in this article.

\vfill\eject
\end{document}